%%Picco december 22 
\newcount\mgnf\newcount\tipi\newcount\tipoformule\newcount\greco

\tipi=2          %uso caratteri: 2=cmcompleti, 1=cmparziali, 0=amparziali
\tipoformule=0   %=0 da numeroparagrafo.numeroformula; se no numero
                 %assoluto

%%%%%%%%%%%%%%%%% EQUAZIONI CON NOMI SIMBOLICI
%%%
%%% Per assegnare un nome simbolico ad una equazione basta
%%% scrivere \Eq(...) o, in \eqalignno, \eq(...) o,
%%% nelle appendici, \Eqa(...) o \eqa(...);
%%% dentro le parentesi e al posto di ... si puo' scrivere qualsiasi commento;
%%% per avere i nomi simbolici segnati a sinistra delle formule si deve
%%% dichiarare il documento come bozza, iniziando il testo con
%%% \BOZZA. Sinonimi: \Eq,\EQ,\EQS; \eq,\eqs; \Eqa,\Eqas;\eqa,\eqas.
%%% All' inizio di ogni paragrafo si devono definire il
%%% numero del paragrafo e della prima formula dichiarando
%%% \numsec=... \numfor=...  (brevetto Eckmannn).
%%% Si possono citare formule seguenti; le corrispondenze fra nomi
%%% simbolici e numeri effettivi sono memorizzate nel file \jobname.aux, che
%%% viene letto all'inizio, se gia' presente. E' possibile citare anche
%%% formule che appaiono in altri file, purche' sia presente il
%%% corrispondente file .aux; basta includere all'inizio l'istruzione
%%%           \include{nomefile}
%%%
%%%%%%%%%%%%%%%%%%%%%%%%%%%%%%%%%%%%%%%%%%%%%%%%%%%%%%%%%%%%%%%

\global\newcount\numsec
\global\newcount\numfor
\global\newcount\numtheo
\global\advance\numtheo by 1

\def\senondefinito#1{\expandafter\ifx\csname#1\endcsname\relax}

\def\SIA #1,#2,#3 {\senondefinito{#1#2}%
\expandafter\xdef\csname #1#2\endcsname{#3}\else
\write16{???? ma #1,#2 e' gia' stato definito !!!!} \fi}

\def\etichetta(#1){(\veroparagrafo.\veraformula)%
\SIA e,#1,(\veroparagrafo.\veraformula) %
\global\advance\numfor by 1%
\write15{\string\FU (#1){\equ(#1)}}%
\write16{ EQ #1 ==> \equ(#1) }}

\def\letichetta(#1){\veroparagrafo.\verotheo
\SIA e,#1,{\veroparagrafo.\verotheo}
\global\advance\numtheo by 1
\write15{\string\FU (#1){\equ(#1)}}
\write16{ Sta \equ(#1) == #1 }}

\def\tetichetta(#1){\veroparagrafo.\veraformula %%%%copy four lines
\SIA e,#1,{(\veroparagrafo.\veraformula)}
\global\advance\numfor by 1
\write15{\string\FU (#1){\equ(#1)}}
\write16{ tag #1 ==> \equ(#1)}}

\def\FU(#1)#2{\SIA fu,#1,#2 }

\def\etichettaa(#1){(A\veroparagrafo.\veraformula)%
\SIA e,#1,(A\veroparagrafo.\veraformula) %
\global\advance\numfor by 1%
\write15{\string\FU (#1){\equ(#1)}}%
\write16{ EQ #1 ==> \equ(#1) }}

\def\BOZZA{
\def\alato(##1){%
 {\rlap{\kern-\hsize\kern-1.4truecm{$\scriptstyle##1$}}}}%
\def\aolado(##1){%
 {%\vtop to \profonditastruttura 
{%\baselineskip
  %\profonditastruttura\vss
 \rlap{\kern-1.4truecm{$\scriptstyle##1$}}}}}
 }

\def\alato(#1){}
\def\aolado(#1){}

\def\veroparagrafo{\number\numsec}
\def\veraformula{\number\numfor}
\def\verotheo{\number\numtheo}

\def\Eq(#1){\eqno{\etichetta(#1)\alato(#1)}}
\def\eq(#1){\etichetta(#1)\alato(#1)}
\def\leq(#1){\leqno{\aolado(#1)\etichetta(#1)}}%%%%%this line for \leqno
\def\teq(#1){\tag{\aolado(#1)\tetichetta(#1)\alato(#1)}}%%%%%this line for\tag
\def\Eqa(#1){\eqno{\etichettaa(#1)\alato(#1)}}
\def\eqa(#1){\etichettaa(#1)\alato(#1)}
\def\eqv(#1){\senondefinito{fu#1}$\clubsuit$#1
\write16{#1 non e' (ancora) definito}%
\else\csname fu#1\endcsname\fi}
\def\equ(#1){\senondefinito{e#1}\eqv(#1)\else\csname e#1\endcsname\fi}

%%%% next six lines by paf (no responsibilities taken)
\def\Lemma(#1){\aolado(#1)Lemma \letichetta(#1)}%
\def\Theorem(#1){{\aolado(#1)Theorem \letichetta(#1)}}%
\def\Proposition(#1){\aolado(#1){Proposition \letichetta(#1)}}%
\def\Corollary(#1){{\aolado(#1)Corollary \letichetta(#1)}}%
\def\Remark(#1){{\noindent\aolado(#1){\bf Remark \letichetta(#1).}}}%
\def\Definition(#1){{\noindent\aolado(#1){\bf Definition 
\letichetta(#1)$\!\!$\hskip-1.6truemm}}}
\def\Example(#1){\aolado(#1) Example \letichetta(#1)$\!\!$\hskip-1.6truemm}

\def\include#1{
\openin13=#1.aux \ifeof13 \relax \else
\input #1.aux \closein13 \fi}

\openin14=\jobname.aux \ifeof14 \relax \else
\input \jobname.aux \closein14 \fi
\openout15=\jobname.aux

%%%%%%%%%%%%%%%%%%%%%%%%%%%%%%%%%%%%%%%%%%%%%%
%%%%%%%%%%%%%%%%%%%%%  Numerazione pagine

{\count255=\time\divide\count255 by 60 \xdef\hourmin{\number\count255}
        \multiply\count255 by-60\advance\count255 by\time
   \xdef\hourmin{\hourmin:\ifnum\count255<10 0\fi\the\count255}}

\def\oramin{\hourmin }

\def\data{\number\day/\ifcase\month\or january \or february \or march \or april
\or may \or june \or july \or august \or september
\or october \or november \or december \fi/\number\year;\ \oramin}

\newcount\pgn \pgn=1
\def\foglio{\number\numsec:\number\pgn
\global\advance\pgn by 1}
\def\foglioa{A\number\numsec:\number\pgn
\global\advance\pgn by 1}

\footline={\rlap{\hbox{\copy200}}\hss\tenrm\folio\hss}
%\footline={\hss\tenrm\folio\hss}

%%%%%%%%
%% am
\def\TIPIO{
\font\setterm=amr7 %\font\settei=ammi7
%\font\settesy=amsy7 \font\settebf=ambx7 %\font\setteit=amit7
%%%%% cambiamenti di formato %%%
\def \settepunti{\def\rm{\fam0\setterm}% passaggio a tipi da 7-punti
\textfont0=\setterm   %\textfont1=\settei
%\textfont2=\settesy   %\textfont3=\setteit
%\textfont\itfam=\setteit  \def\it{\fam\itfam\setteit}
%\textfont\bffam=\settebf  \def\bf{\fam\bffam\settebf}
\normalbaselineskip=9pt\normalbaselines\rm }\let\nota=\settepunti}
%%%%%%%

%%cm completo
\def\TIPITOT{
\font\twelverm=cmr12
\font\twelvei=cmmi12
\font\twelvesy=cmsy10 scaled\magstep1
\font\twelveex=cmex10 scaled\magstep1
\font\twelveit=cmti12
\font\twelvett=cmtt12
\font\twelvebf=cmbx12
\font\twelvesl=cmsl12
\font\ninerm=cmr9
\font\ninesy=cmsy9
\font\eightrm=cmr8
\font\eighti=cmmi8
\font\eightsy=cmsy8
\font\eightbf=cmbx8
\font\eighttt=cmtt8
\font\eightsl=cmsl8
\font\eightit=cmti8
\font\sixrm=cmr6
\font\sixbf=cmbx6
\font\sixi=cmmi6
\font\sixsy=cmsy6
%%%%%%%%%%%%%%%%%%%%%%%%%%%%%%%%%%%%%%%
\font\twelvetruecmr=cmr10 scaled\magstep1
\font\twelvetruecmsy=cmsy10 scaled\magstep1
\font\tentruecmr=cmr10
\font\tentruecmsy=cmsy10
\font\eighttruecmr=cmr8
\font\eighttruecmsy=cmsy8
\font\seventruecmr=cmr7
\font\seventruecmsy=cmsy7
\font\sixtruecmr=cmr6
\font\sixtruecmsy=cmsy6
\font\fivetruecmr=cmr5
\font\fivetruecmsy=cmsy5
%%%% definizioni per 10pt %%%%%%%%
\textfont\truecmr=\tentruecmr
\scriptfont\truecmr=\seventruecmr
\scriptscriptfont\truecmr=\fivetruecmr
\textfont\truecmsy=\tentruecmsy
\scriptfont\truecmsy=\seventruecmsy
\scriptscriptfont\truecmr=\fivetruecmr
\scriptscriptfont\truecmsy=\fivetruecmsy
%%%%% cambio grandezza %%%%%%
\def \eightpoint{\def\rm{\fam0\eightrm}% switch to 8-point type
\textfont0=\eightrm \scriptfont0=\sixrm \scriptscriptfont0=\fiverm
\textfont1=\eighti \scriptfont1=\sixi   \scriptscriptfont1=\fivei
\textfont2=\eightsy \scriptfont2=\sixsy   \scriptscriptfont2=\fivesy
\textfont3=\tenex \scriptfont3=\tenex   \scriptscriptfont3=\tenex
\textfont\itfam=\eightit  \def\it{\fam\itfam\eightit}%
\textfont\slfam=\eightsl  \def\sl{\fam\slfam\eightsl}%
\textfont\ttfam=\eighttt  \def\tt{\fam\ttfam\eighttt}%
\textfont\bffam=\eightbf  \scriptfont\bffam=\sixbf
\scriptscriptfont\bffam=\fivebf  \def\bf{\fam\bffam\eightbf}%
\tt \ttglue=.5em plus.25em minus.15em
\setbox\strutbox=\hbox{\vrule height7pt depth2pt width0pt}%
\normalbaselineskip=9pt
\let\sc=\sixrm  \let\big=\eightbig  \normalbaselines\rm
\textfont\truecmr=\eighttruecmr
\scriptfont\truecmr=\sixtruecmr
\scriptscriptfont\truecmr=\fivetruecmr
\textfont\truecmsy=\eighttruecmsy
\scriptfont\truecmsy=\sixtruecmsy }\let\nota=\eightpoint}

\newfam\msbfam   %per uso in \TIPITOT
\newfam\truecmr  %per uso in \TIPITOT
\newfam\truecmsy %per uso in \TIPITOT
%%%%%%%%%%%%%%%%%%%%%%%%%%%%%%%
%%Scelta dei caratteri
%\newcount\tipi \tipi=0   %e' definito all'inizio
\newskip\ttglue
\ifnum\tipi=0\TIPIO \else\ifnum\tipi=1 \TIPI\else \TIPITOT\fi\fi

\def\a{\alpha}
\def\b{\beta}

\def\g{\gamma}

\def\l{\lambda}
\def\r{\rho}
\def\s{\sigma}
\def\t{\tau}
\def\th{\theta}

\def\z{\zeta}
\def\o{\omega}
\def\D{\Delta}
\def\L{\Lambda}
\def\G{\Gamma}
\def\O{\Omega}
\def\S{\Sigma}

\def\E{{I\kern-.25em{E}}}
\def\N{{I\kern-.25em{N}}}
\def\M{{I\kern-.25em{M}}}
\def\R{{I\kern-.25em{R}}}
\def\Z{{Z\kern-.425em{Z}}}
\def\1{{1\kern-.25em\hbox{\rm I}}}
\def\eu{{1\kern-.25em\hbox{\sm I}}}

\def\C{{I\kern-.64em{C}}}
\def\P{{I\kern-.25em{P}}}
\def\eop{{ \vrule height7pt width7pt depth0pt}\par\bigskip}

%\def\P{\hskip.2em\hbox{\rm P\kern-0.8em{I}\hskip.7em}}

% Spezielle Definitionen

\def\AA{{\cal A}}
\def\BB{{\cal B}}

\def\SS{{\cal S}}

\def\XX{{\cal X}}

\def\RR{{\cal R}}

\def\chap #1#2{\line{\ch #1\hfill}\numsec=#2\numfor=1}

\def\sqr#1#2{{\vcenter{\vbox{\hrule height.#2pt
     \hbox{\vrule width.#2pt height#1pt \kern#1pt
   \vrule width.#2pt}\hrule height.#2pt}}}}

%   Non-character macros

\newcount\foot
\foot=1
\def\note#1{\footnote{${}^{\number\foot}$}{\ftn #1}\advance\foot by 1}
\def\tag #1{\eqno{\hbox{\rm(#1)}}}
\def\frac#1#2{{#1\over #2}}

\def\text#1{\quad{\hbox{#1}}\quad}

\def\proof{{\noindent\pr Proof: }}

\def\thanks{\noindent{\bf Aknowledgements: }}
\font\pr=cmbxsl10
%\font\thbf=cmcsc10 scaled\magstep1

% Font-Definitions

\font\ch=cmbx12
\font\ftn=cmr8

\font\it=cmti10
\font\bf=cmbx10
\font\sm=cmr7

%%%%%%%%References macros start here%%%%%%%%%%
%
\catcode`\X=12\catcode`\@=11
\def\n@wcount{\alloc@0\count\countdef\insc@unt}
\def\n@wwrite{\alloc@7\write\chardef\sixt@@n}
\def\n@wread{\alloc@6\read\chardef\sixt@@n}
\def\crossrefs#1{\ifx\alltgs#1\let\tr@ce=\alltgs\else\def\tr@ce{#1,}\fi
   \n@wwrite\cit@tionsout\openout\cit@tionsout=\jobname.cit 
   \write\cit@tionsout{\tr@ce}\expandafter\setfl@gs\tr@ce,}
\def\setfl@gs#1,{\def\@{#1}\ifx\@\empty\let\next=\relax
   \else\let\next=\setfl@gs\expandafter\xdef
   \csname#1tr@cetrue\endcsname{}\fi\next}
\newcount\sectno\sectno=0\newcount\subsectno\subsectno=0\def\r@s@t{\relax}
\def\resetall{\global\advance\sectno by 1\subsectno=0
  \gdef\firstpart{\number\sectno}\r@s@t}
\def\resetsub{\global\advance\subsectno by 1
   \gdef\firstpart{\number\sectno.\number\subsectno}\r@s@t}
\def\v@idline{\par}\def\firstpart{\number\sectno}
\def\l@c@l#1X{\firstpart.#1}\def\gl@b@l#1X{#1}\def\t@d@l#1X{{}}
\def\m@ketag#1#2{\expandafter\n@wcount\csname#2tagno\endcsname
     \csname#2tagno\endcsname=0\let\tail=\alltgs\xdef\alltgs{\tail#2,}%
  \ifx#1\l@c@l\let\tail=\r@s@t\xdef\r@s@t{\csname#2tagno\endcsname=0\tail}\fi
   \expandafter\gdef\csname#2cite\endcsname##1{\expandafter
 %the following line was replaced by the subseqent one, DNA 7/6/89
  %  \ifx\csname#2tag##1\endcsname\relax?\else\csname#2tag##1\endcsname\fi
     \ifx\csname#2tag##1\endcsname\relax?\else{\rm\csname#2tag##1\endcsname}\fi
    \expandafter\ifx\csname#2tr@cetrue\endcsname\relax\else
     \write\cit@tionsout{#2tag ##1 cited on page \folio.}\fi}%
   \expandafter\gdef\csname#2page\endcsname##1{\expandafter
     \ifx\csname#2page##1\endcsname\relax?\else\csname#2page##1\endcsname\fi
     \expandafter\ifx\csname#2tr@cetrue\endcsname\relax\else
     \write\cit@tionsout{#2tag ##1 cited on page \folio.}\fi}%
   \expandafter\gdef\csname#2tag\endcsname##1{\global\advance
     \csname#2tagno\endcsname by 1%
   \expandafter\ifx\csname#2check##1\endcsname\relax\else%
\fi%      \immediate\write16{Warning: #2tag ##1 used more than once.}\fi
   \expandafter\xdef\csname#2check##1\endcsname{}%
   \expandafter\xdef\csname#2tag##1\endcsname
     {#1\number\csname#2tagno\endcsnameX}%
   \write\t@gsout{#2tag ##1 assigned number \csname#2tag##1\endcsname\space
      on page \number\count0.}%
   \csname#2tag##1\endcsname}}%
\def\m@kecs #1tag #2 assigned number #3 on page #4.%
   {\expandafter\gdef\csname#1tag#2\endcsname{#3}
   \expandafter\gdef\csname#1page#2\endcsname{#4}}
\def\re@der{\ifeof\t@gsin\let\next=\relax\else
    \read\t@gsin to\t@gline\ifx\t@gline\v@idline\else
    \expandafter\m@kecs \t@gline\fi\let \next=\re@der\fi\next}
\def\t@gs#1{\def\alltgs{}\m@ketag#1e\m@ketag#1s\m@ketag\t@d@l p
    \m@ketag\gl@b@l r \n@wread\t@gsin\openin\t@gsin=\jobname.tgs \re@der
    \closein\t@gsin\n@wwrite\t@gsout\openout\t@gsout=\jobname.tgs }
\outer\def\localtags{\t@gs\l@c@l}
\outer\def\globaltags{\t@gs\gl@b@l}
\outer\def\newlocaltag#1{\m@ketag\l@c@l{#1}}
\outer\def\newglobaltag#1{\m@ketag\gl@b@l{#1}}

\def\t@gsoff#1,{\def\@{#1}\ifx\@\empty\let\next=\relax\else\let\next=\t@gsoff
   \expandafter\gdef\csname#1cite\endcsname{\relax}
   \expandafter\gdef\csname#1page\endcsname##1{?}
   \expandafter\gdef\csname#1tag\endcsname{\relax}\fi\next}
\def\verbatimtags{\let\ift@gs=\iffalse\ifx\alltgs\relax\else
   \expandafter\t@gsoff\alltgs,\fi}
\catcode`\X=11 \catcode`\@=\active
\localtags
%%%%%%%%%%%%references macro end here%%%%%%%%%%%%%%%
%
%%%%%%%%%%%%%%%%%end of macros%%%%%%%%%%%%%%%%%%%%%%
%\BOZZA
\setbox200\hbox{$\scriptscriptstyle \data $}
\global\newcount\numpunt
%\magnification=\magstephalf
%\magnification=\magstep1
\hoffset=0.cm
%\voffset=-.2truecm
%\hoffset=-.2truecm
%\hsize=16.5truecm 
%\vsize=24truecm
\baselineskip=14pt  
\parindent=12pt
\lineskip=4pt\lineskiplimit=0.1pt
\parskip=0.1pt plus1pt

\hyphenation{small}

%\input newdef.tex
%  \BOZZA

\catcode`\@=11

\centerline {\bf  Phase Transition in the 1d Random Field   Ising Model  with  long range interaction.  
 \footnote* 
   {\eightrm Supported by:    GDRE 224 GREFI-MEFI, CNRS-INdAM. P.P was also partially supported by INdAM
program Professori Visitatori 2007; M.C and E.O were partially supported by Prin07: 20078XYHYS.    
}}
\vskip.5cm
 \centerline{ 
Marzio Cassandro \footnote{$^1$}{\eightrm 
Dipartimento di Fisica,  
Universit\'a di Roma ``La Sapienza'',
P.le A. Moro, 00185 Roma, Italy. cassandro@roma1.infn.it}
\hskip.2cm 
Enza Orlandi \footnote{$^2$}{\eightrm 
Dipartimento di Matematica, Universit\'a di Roma Tre, L.go S.Murialdo 1,
00146 Roma, Italy. orlandi@mat.uniroma3.it} and 
Pierre Picco \footnote{$^3$}{\eightrm LATP, CMI, UMR 6632,  CNRS,
Universit\'e de Provence,  39 rue Frederic  Joliot Curie,  13453
Marseille Cedex 13, France. picco@cmi.univ-mrs.fr }  }

 \footnote{}{\eightrm {\eightit AMS 2000 Mathematics Subject Classification}:
Primary 60K35, secondary 82B20,82B43.}
\footnote{}{\eightrm {\eightit Key Words}: phase transition, long--range interaction, random  field. 
 }

\vskip .5cm
{ \bf Abstract}
We  study one--dimensional Ising spin systems with ferromagnetic, long--range interaction   
decaying as $n^{-2+\a}$, $\a \in (\frac 12,  \frac {\ln 3}  {\ln2}-1)$,  in the presence of external random fields.  
 We assume that  the random fields are  given by a collection  of  symmetric, independent,  identically distributed 
real  random variables, gaussian or subgaussian. 
   We show, for temperature  and    
 strength of the randomness (variance)  small enough,  with $\P=1$
with respect to  the 
%{distribution} of the
 random fields, 
that  there are at least two distinct   extremal Gibbs measures. 
\bigskip \bigskip

\chap {1 Introduction}1
\numsec= 1 \numfor= 1

%{\bf  It is well known that in order to have phase transition for one dimensional bounded  %spin systems one needs an interaction  that is  ``so long range'' that it has divergent %first moment.  On the other  hand  necessary condition for the existence of the %thermodynamical   potential  requires the interaction  to  be summable  see %[\rcite{GM}].   }

It is well known that the one dimensional  ferromagnetic Ising model exhibits  a phase transition when the forces are sufficiently long range.  A fundamental   work on the subject is due 
to Dyson [\rcite{Dy}].   He    proved,  by comparison  to a    hierarchical model,      that  for   a two body  interaction   $J(n)=\frac {\ln \ln (n+3)}  {n^{2} +1}$,  where  $n$ denotes the distance,    there is spontaneous magnetization at low enough  temperature.

On the other hand 
Rogers \&Thompson [\rcite{RT}] proved that the spontaneous magnetization vanishes for all temperatures  when
$$ \lim_{N \to \infty} \frac 1 {[\ln N]^{\frac 12}} \sum_{n=1}^N  n J(n) =0. $$ 
Later,    Fr\"ohlich \& Spencer [\rcite{FS}] 
  proved   the existence of spontaneous magnetization  when $J(n)=n^{-2}$.
For the same model
 Aizenman,  Chayes,  Chayes, \&  Newman [\rcite{ACCN}]   proved
 %, always in the case $J(n)=n^{-2}$, 
  the 
 discontinuity   of the magnetization at the critical temperature, the so-called Thouless effect. 
   %It is well known that for one dimensional bounded  spin system with long range %interactions
%one has existence of thermodynamics, see [\rcite{GM}],  when 
%the interaction of any  single site  with the rest of the spins in lattice $\Z$ is finite. 
%For two body interactions which decays as  $|i-j|^{-2+\a}$ % this means that $\a<1$. 
When  $J(n)=n^{-2+\a}$, 
 $\a<0$ there is only one Gibbs state  [\rcite{Ru},\rcite{D0},\rcite{D1}] and  
the free energy is analytic in the thermodynamic parameters, see [\rcite{D2}].
     More recently the notion of contours  introduced in    [\rcite {FS}] was implemented in  [\rcite {CFMP}], 
by giving a graphical description of the spin configurations
% more explicit and 
better suited for further
generalizations.  The case studied in [\rcite{CFMP}] covers the regime $0\le \a \le  (\ln 3/\ln 2)-1$.
By  applying   Griffiths inequalities    the existence of a phase transition  
in the full interval  $0\le \a<1$ can be deduced either by  [\rcite{CFMP}], or by   [\rcite{FS}].

A natural extension of this analysis is its application to  disordered systems.  
One of the simplest  prototype models for disordered spin systems is 
obtained by  adding  random magnetic fields, say gaussian independent 
identically distributed with mean zero and finite variance. The problem  
of  (lower) critical dimension for  the    d--dimensional  Random Field Ising Model 
was very challenging 
at the end of the eighties since the physical literature predicted conflicting results.  
 For  %the  short 
   finite   range interaction the problem was  rigorously solved by two complementary articles,  
 Bricmont \& Kupiainen [\rcite{BK}] and Aizenman \& Wehr
[\rcite{AW}]. In [\rcite{BK}]  a  
renormalization group argument was used to show that if 
$d\ge 3$ and the variance of the random magnetic field is small enough then almost surely 
there are  at least two distinct Gibbs states (the plus and the minus Gibbs states).
In [\rcite{AW}] it was proved that for $d\le 2$, almost surely there is  an unique Gibbs state. 
  The guide lines of these proofs are suggested by a   heuristic argument  due to  Imry \& Ma [\rcite{IM}]. 

 In the long--range one--dimensional setting the Imry \& Ma argument is the following: the deterministic
cost to create a run of $-1$ in an interval of length $L$ with respect to the state
 made of $+1$ at each site, 
is of order $L^{\a}$,  while the cumulative effect of the random field inside this interval 
 is just $L^{1/2}$.
So when $0\le \a \le1/2$ the randomness is dominant  and there is no phase transition. 
This 
has been   proved  by Aizenman \& Wehr   [\rcite {AW}].  They  show   that   
the Gibbs state is unique for almost  all realizations of the randomness. 

When $1/2<\a<1$,    the above Imry \& Ma argument suggests  the existence of a phase transition 
since the deterministic part is dominant with respect to the random part  
as in the case of the three--dimensional random field Ising model.
However a rigorous result is this direction was missing.

In this paper
we study the random field   one--dimensional   Ising model 
with long range  interaction    $n^{-2+\a}$, $   \a\in (\frac 12,  \frac { \ln 3}{\ln 2}-1)\simeq (\frac 12,  \frac {58}  {100})$.   
We assume that the random field  $ h[\o]: = \{ h_i [\o], i \in \Z \} $ is given 
by a collection  of independent random variables, with mean zero and symmetrically  distributed. 
We take $h_i [\o]=\pm 1$ with $p=\frac 12$  and we introduce the strength  
parameter $\th$.  However one could take different distributions,  for example gaussian  
distribution with  mean zero and variance $\theta^2$, or subgaussian. In fact  all that  is needed 
is   $E(e^{t h_1})\le e^{c\th^2 t^2}$ for some positive constant $c$ and for all $t\in \R$. 
We prove   that for  $\frac 12 <\a <  \frac { \ln 3}{\ln 2}-1 $ the situation is analogous to the three-dimensional
short range random field Ising model:  for temperature and variance of the randomness small enough,  
with $\P=1$ with respect to the randomness, 
there exist at least two distinct infinite volume Gibbs states, namely  the
$ \mu^+ [\o]$ and the  $ \mu^- [\o]$ Gibbs states.
The proof is  based on the representation of the system in term of the     contours as defined 
in [\rcite{CFMP}]. A  Peierls argument is obtained  by using  
the lower bound of the deterministic part of the cost to erase a contour 
%as 
%done  in [\rcite{CFMP}] 
and   controlling   the contribution of the stochastic part.
This control is done  applying  an exponential Markov inequality and 
the so-called Yurinski's martingale difference sequences method.
  We do not need to use  any coarse-grained contours as in [\rcite{BK}], a fact that simplifies the proof. 
%We do not use any coarse-grained contours as in [\rcite{BK}], a fact that simplifies the % proof. 
 %This  is possible in our case   because 
In the one dimensional case the  
 contours   can be described in terms of intervals    and  the Imry\& Ma
argument  can be  implemented.  Namely in our case   bad configurations of the random magnetic field, the ones for which  the naive Imry \& Ma argument fails, are treated probabilistically.   A kind of  energy entropy argument is  successfully used, see \eqv (En.1), to prove that they can be neglected.   In 3 dimensions this  specific energy entropy   argument
% to treat these bad  configurations of the magnetic fields
  fails. 
  The coarse grained contours in [\rcite{BK}] allow  to control these bad contours
on various length scales by using a renormalization group argument. As a by-product an estimate on the  decay of the truncated two point correlation functions is given in [\rcite{BK}].  Our method  does  not give any information on this decay. 
%Therefore we do not think that our method can be used to give an alternative proof of %the existence of a phase transition in three dimension [BK]..
Therefore we  do not think that
% our method 
it can be directly applied to give an alternative proof of Bricmont \& Kupiainen  results  [\rcite{BK}].    
% In higher dimensions the contours have  a more complex  topological structure
%and it is  convenient  to introduce a coarse graining procedure  (see  [\rcite{BK}]) to
%single out  and control the two conflicting contributions:
% stochastic and deterministic 
   For $ \a \in [  (\ln 3/\ln 2)-1, 1) $ we still expect the same result  to hold but we are not able to prove it.   In this case 
  the lower  bound for the deterministic contribution to the cost of erasing a contour  does not hold, see  Lemma  \eqv (CFMP1).  
  Known correlations inequalities  are not relevant  to  treat  this range of values of $\a$  as in the case where the random field is absent. 
   %On the other hand to prove  our result  for    $ \a \in [  (\ln 3/\ln 2)-1, 1) $ would %require an  extension of the estimates proven in [\rcite{CFMP}]  for the above values of $\a$.  This would require new ideas. 
  % The  extension   for    $ \a \in [  (\ln 3/\ln 2)-1, 1) $ of the estimates proven  in  [\rcite{CFMP}] which we use  to lower bound  the deterministic   cost to erase a contour,   requires  the use of some %new methods.  %  The main difficulty we face is that we are not able   to  extend   the deterministic %estimates of the cost of the contours, see  Theorem \eqv   (2CFMP),  which  holds   %only
  %for  $ 0 <\a < (\ln 3/\ln 2)-1$.      
       \medskip 
 
\vskip 1truecm
\chap{2 Model, notations  and main results }2
\numsec= 2
\numfor= 1
\numtheo=1 

\medskip
\noindent{\bf 2.1. The model and the main results }
\medskip
Let $(\O,\BB,\P)$ be a probability space on which we  define 
$h \equiv \{h_i\}_{i\in \Z}$, a family of
independent, identically distributed Bernoulli random variables with
$ \P[h_i=+1]=\P[h_i=-1]=1/2$.
The spin configuration space is $\SS\equiv  \{-1,+1\}^\Z$.
If $\s \in \SS$ and $i\in \Z$,
$\s_i$ represents the value of the spin at site $i$. 
The pair interaction among spins is given by  $J(|i-j|)$ defined as 
following\footnote{$^1$}{\eightrm  The condition $J(1) >>1  $ is essential
 to apply the results of   [\rcite {CFMP}], reported in Subsection 2.2.}:
 $$J(n)= \left \{ \eqalign {& J(1) >>1      \cr &
 \frac 1 {n^{2- \alpha} } \quad \hbox {if}\quad  n >1, \quad \hbox {with} \quad \alpha \in [0,1).} \right. $$ 
 For $\L \subseteq \Z$ we set $\SS_\L=\{-1,+1\}^\L$; its elements
are   denoted by $\s_\L$; also,
if $\s \in \SS$, $\s_\L$ denotes its restriction to $\L$. Given
$\L\subset \Z$  finite and a realization of the magnetic fields,
the   Hamiltonian in the  volume $\L$, with $\t=\pm 1$  boundary conditions,
is the  random variable on $(\O,\AA,\P)$
given by
$$
H^{\t} (\s_\L)[\o]= H^\t_0(\s_\L) +\th G(\s_\L)[\o]
\Eq(2.1)
$$
where
$$ 
H^\t_0(\sigma_\Lambda):= \frac 12 \sum_{(i,j) \in \L \times \L}
 J(|i-j|)(1- \s_i \s_j)   +
\sum_{i\in \L} \sum_{j\in \L^c} J (|i-j|)
(1-\t \s_i ), 
\Eq(2.1a)
 $$
and 
$$
G(\sigma_\Lambda)[\o]:=- \sum_{i\in \L} h_i[\o] \s_i. 
\Eq(2.1ab)
$$

In the following we drop
the $\o$ from the notation. 
The corresponding {\sl Gibbs measure} on the finite volume $\L$,
at inverse temperature $\b>0$ and  $+$
boundary condition  is then a random variable with value 
on the space of probability measures on $\SS_\L$
defined by
$$
\mu^+_{ \L}(\s_\L)
= \frac 1{Z^+_{ \L}} \exp\{-\b H^+(\s_\L)\} \quad \quad
\s_\L \in \SS_\L,
\Eq(2.3)
$$
where  $Z^+_{\L}$ is the normalization factor. 
Using  FKG inequalities, 
one can construct  with $ \P=1$ 
the  infinite volume Gibbs measure  $ \mu^+[\o]$ as limits of local specifications with 
 homogeneous  plus boundary conditions along any deterministic sequence of increasing and 
absorbing finite volumes $\L_n$.  Of course the same   holds with minus boundary conditions,
see for example  Theorem 7.2.2 in  [\rcite {Bo}]  or Theorem IV.6.5 in 
[\rcite{E}].    
 The main results are the following. 
 \vskip0.4cm \noindent  
{\bf \Theorem(1)  } {\it     Let $ \alpha \in (\frac 12, \frac {\ln3} {\ln2} -1)$
and  $$\zeta=\zeta(\a)=1-2(2^\a-1)>0.
\Eq(fo.2bis)$$
There exist  positive   
 $\theta_0:= \theta_0(\alpha)>0 $ and  $ \b_0:=\b_0(\a)>0 $  so that for  
 $0<\theta \le \theta_0 $  and     $\b \ge  \b_0 $      there exists   $\O_1 \subset \O $ such that
$$
\P[\O_1] \ge 
1-   e^{-\frac {\bar b} {200}  },     
\Eq(5.02)$$
 and for any  $\o\in \O_1$,  $$
\mu^+ \Big(\{ \s_0=-1 \}\Big)[\o]
<  e^{-\frac {\bar b} {200}  }  \Eq (R1)  
$$
where  
$$ \bar b = \min(\frac {\beta \zeta} 4,  \frac {\z^2} {2^{10}\theta^2}).  \Eq (SS5) $$
} 
\vskip0.5cm \noindent 
 %\noindent {\bf Remark:} If we prove phase transition for
% $ \alpha \in (\frac 12, \frac {\ln3} {\ln2} -1)$ then FKG inequalities 
%imply  that the  result holds  for $ \alpha \in (\frac 12, 1)$.  
%\medskip
\noindent {\bf Remark:} 
 Since the translation invariant, $ \BB$ measurable event 
$A\equiv \{\exists i \in \Z:  \mu^+[\o] (\s_i=+1) > 1-e^{-\frac{\bar b}{200}} \}$
has strictly positive probability, see \eqv(5.02) and \eqv(R1),  by ergodicity $\P[A]=1$.
Therefore almost surely the two extremal Gibbs states $\mu^{\pm}[\o]$ are distinct.
 \medskip
 \noindent 
 The proof of Theorem  \eqv (1)   is given in Section 3. In the next subsection we recall the definition of contours 
and in Section 4 we prove the main probabilistic estimate.
 \vskip0.5cm
\medskip
\noindent{\bf 2.2.  Geometrical description of the spin configurations}
\medskip

We will follow  the geometrical description of the spin configuration  presented in  [\rcite {CFMP}]
and  use the same notations.
We will consider homogeneous boundary conditions, i.e the spins in the
boundary conditions are either all $+1$ or all $-1$. Actually we will restrict ourself to
$+$ boundary conditions and  consider spin   configurations $ \s= \{ \s_i, i\in \Z \}\in \XX_+ $ 
so that $\s_i=+1$ for all $|i|$ large enough.

In one dimension an interface at $(x, x+1)$ means  $\s_{x} \s_{x+1}= -1$.   
Due to the above choice of the boundary conditions, any $\s \in \XX_+ $ has a finite, 
even number of  interfaces.   The precise location of the interface is immaterial and 
 this  fact has been  used  to choose the interface points as  follows:   For
 all $x \in \Z$ so that $(x,x+1)$ is an interface take the  location of the 
interface  to be a point  inside the interval  
$[x+\frac 12- \frac 1 {100}, x+\frac 12+ \frac 1 {100} ] $, with the property that for 
any four distinct points $r_i$,  $i=1, \dots, 4$ $|r_1-r_2| \neq |r_3-r_4|$.   This choice 
is done once for all so that the interface between $x$ and $x+1$ is uniquely fixed.     
 Draw  from each one of  these interfaces points  two 
lines  forming  respectively   an angle of $ \frac \pi 4$  and of   $ \frac 3 4 \pi  $ with 
the $\Z$ line.    We have thus a bunch of growing $\vee-$ lines each one emanating from an
 interface point.    Once two $\vee-$ lines meet, they are frozen and stop their growth. 
 The other two lines emanating  from the the same interface points are {erased}.
The $\vee-$ lines 
 emanating from others points  keep growing.  The collision of the two lines is represented
 graphically by a triangle whose basis is the line joining the two interfaces points and whose 
sides are the two segment of the $\vee-$ lines which meet.   The choice done of the location 
of the interface points ensure that collisions occur one at a time so that the above  definition 
is unambiguous.   In general  there might be    triangles inside triangles.   
  The endpoints of the triangles 
  are    suitable coupled    pairs of interfaces points.    The  graphical representation  
just described   maps  each spin configuration in $ \XX_+ $ to a set of triangles.

  \vskip0.5cm
  \noindent {\bf Notation }
 {\it  Triangles will be usually denoted by $T$, the collection of triangles constructed as above 
  by $ \{ \underline T \}$ and we will write 
  $$ |T| = \hbox {cardinality of} \quad T \cap \Z =  \hbox {mass of} \quad T, $$
 and by   $\hbox {supp} (T) \subset  \R $ the   basis of the triangle.}
 
 \medskip \noindent
  We have thus represented a configuration $\s \in \XX_+ $ as a collection of 
$\underline T = (T_1, \dots, T_n)$.  The above construction defines a one to
 one map from $ \XX_+ $ onto $ \{ \underline T \}$.  It is easy to see that a triangle 
configuration $\underline T$ belongs to $ \{ \underline T \}$ iff for any pair $T$ and 
$T'$ in $\underline T$ 
 $$ {\rm dist} (T,T')\ge \min\{ |T|, |T'|\} .\Eq (Ma1)$$ 
We say that  two collections of triangles $\underline S'$  and $\underline S$   are
 compatible  and we denote it by  $\underline S' \simeq \underline S$  iff 
$ \underline S' \cup \underline S \in  \{ \underline T \}$ 
({\it i.e.} 
there exists a configuration in $\XX_+$ such that its corresponding collection  of triangles 
is the collection made of all triangles that are in $\underline S'$ or in $\underline S$.)
By an abuse of notation, we write 
   $$  H^+_0(\underline T)= H^+_0(\sigma),  \quad G(\sigma(\underline T) )[\o]= G(\sigma)[\o], 
 \quad \s \in \XX_+ \iff \underline T \in   \{ \underline T \}. $$ 
     \vskip0.5cm

  \noindent {\bf \Definition  (3)}   { \bf The energy difference } {\it 
 Given  two compatible   collections   of triangles  $\underline S \simeq   \underline T$,    we denote
     $$
H^+(\underline S| \underline T):=   H^+(\underline S\cup \underline T)-H^+( \underline T). \Eq (D3) $$}

Let  $ \underline T= (T_1, \dots, T_n)$ with $|T_i| \le |T_{i+1}|$ then using \eqv(D3) one has
$$
H^+(\underline T)=H^+(T_1|\underline T\setminus T_1)+ H^+(\underline T\setminus T_1).
\Eq(D3bis)
$$

The following Lemma proved  in  [\rcite {CFMP}], see  Lemma 2.1 there,  gives a lower bound on the cost to ``erase'' 
triangles sequentially starting from the smallest ones.

\medskip \noindent 
 \noindent  { \bf \Lemma (CFMP1)   [\rcite {CFMP}] } {\it For $ \alpha \in (0, \frac {\ln3} {\ln2} -1)$   
and  $ \zeta:= \zeta(\a)$ as  defined  in \eqv (fo.2bis)
one has 
$$ H^+_0(T_1| \underline T \setminus  T_1 )  \ge \zeta |T_1|^\a,  \Eq (Ma2) $$ 
and by iteration, 
for  
 any $1\le i \le n$}
   $$ H^+_0(\cup_{\ell=1}^i T_{\ell} | \underline T \setminus [\cup_{\ell=1}^i T_{\ell} ]) 
 \ge \zeta \sum_{\ell=1}^i |T_\ell|^\a. \Eq (Ma2a) $$ 
The  estimate \eqv (Ma2a)  involves contributions coming from the full set of triangles 
associated to a given spin configuration, starting from the triangle having the smallest 
mass.  To  implement a Peierls bound in our set up we need to ``localize'' the estimates 
 to compute the weight  of a  triangle or of a finite set of triangles in a generic configuration.  
In order to   do this  [\rcite {CFMP}]       introduced  the notion of contours as clusters of
 nearby triangles sufficiently far away from all other triangles.
 %% i.e weakly coupled to the surrounding landscape of triangles.    
  \medskip \noindent {\bf Contours} 
   A contour $ \G$ is a collection $\underline T$ of  triangles related by a 
hierarchical network of connections  controlled by a positive number $C$, see \eqv (SS1), under which all the triangles of a contour become mutually 
connected.     We denote by $T(\G)$ the triangle  whose basis is the smallest interval which contains all the 
triangles of the contour.  The right and left endpoints of $T(\G) \cap \Z$ are denoted by $x_{\pm} (\G)$.   
 We denote $|\G|$ the  mass of the contour  $\G$
 $$ |\G|= \sum_{T \in \G}|T| $$
 i.e.  $ |\G|$ is   the sum of the masses of all the triangles  belonging to $\G$.  
We denote  by   $ \RR (\cdot)$  the algorithm 
%In  [\rcite {CFMP}]   it has been defined   an algorithm $ \RR (\underline T)$ on %$\{\underline T\}$ 
which 
associates to any configuration $ \underline T $  a configuration $ \{ \G_j\}$ of contours with the following 
properties.   
   \vskip0.5cm   \noindent
   { \bf P.0}  {\it   Let $ \RR (\underline T) = ( \G_1, \dots, \G_n)$, $ \G_i= \{ T_{j,i}, 1 \le j\le k_i\}$, 
then $\underline T=  \{ T_{j,i}, 1\le i \le n, 1 \le j\le k_i\}$} 
    
   \vskip0.5cm  \noindent
   { \bf P.1}  {\it Contours are well separated from each other.}  Any pair $ \G \neq \G'$ 
 verifies one of the following alternatives.
  
  $$       T (\G) \cap T (\G') = \emptyset $$
  i.e.  $[x_{-} (\G), x_{+} (\G)] \cap  [x_{-} (\G'), x_{+} (\G')] = \emptyset$,  in which case 
  $$ dist (\G, \G'):= \min_{T \in \G, T' \in \G'} dist ( T,T') > C \left \{ |\G|^3, |\G'|^3\right \} \Eq (SS1)  $$
  where $C$ is a positive number.% and it is chosen so that  \eqv (SS2) holds.  
 If 
   $$     T (\G) \cap T (\G')\neq  \emptyset, $$
   then either $T(\G)\subset T(\G')$ or $T(\G')\subset T(\G)$; moreover, supposing for instance
 that the former case is verified, (in which case we call $\G$ an inner contour) then for any 
triangle $ T'_i \in \G'$, either $T(\G)\subset T'_i$ or  $T(\G)\cap  T'_i= \emptyset $  and 
   $$ dist (\G, \G') >C |\G|^3, \quad \hbox {if} \quad  T(\G)\subset T(\G').  \Eq (SS3)$$
   \vskip0.5cm
   \noindent  
   { \bf P.2}  {\it Independence.}  Let $\{ \underline T^{(1)}, \dots, \underline T^{(k)}\}$,
 be $k>1$ configurations of triangles;    $\RR ( \underline T^{(i)}) = \{ \G_j^{(i)}, j=1,\dots, n_i\}$
 the contours of the configurations $\underline T^{(i)}$. Then if any distinct $ \G_j^{(i)}$ and 
 $\G_{j'}^{(i')}$ satisfies {\bf P.1}, 
   $$\RR ( \underline T^{(1)}, \dots, \underline T^{(k)} )= \{  \G_j^{(i)}, j=1,\dots, n_i; i=1, \dots, k \}. $$
  %In  [\rcite {CFMP}]  it has been proved that not  only  {\bf P.0},  {\bf P.1} and   {\bf P.2} can be actually
% implemented by some but also that such  algorithm is unique   
  As proven in   [\rcite {CFMP}], the  algorithm $ \RR (\cdot)$ having properties  {\bf P.0},  {\bf P.1} and   {\bf P.2}   is unique  and therefore there is a 
bijection between  families of triangles  and contours. 
   Next we report the   estimates proven in   [\rcite {CFMP}] which are essential  for  this paper. 
 \vskip0.5cm
  \noindent  { \bf \Theorem (1CFMP) [\rcite {CFMP}] } {\it  Let 
 $\alpha \in (0, \frac {\ln3} {\ln2} -1)$ and the constant $C$ in the 
definition of the contours, see \eqv (SS1),    be so large that 
  $$ \sum_{m\ge 1} \frac {4m} {[C m]^3} \le \frac 12,\Eq (SS2)$$
  where  $ [x ]$ denotes  the integer part of $x$. 
For any $ \underline  T \in  \{ \underline T\}$,  let $\G_0\in \RR (\underline T)$ 
be a contour, $\underline S^{(0)}$ the triangles in $\G_0$ and     $  \zeta(\a)$ as in  \eqv  (fo.2bis) 
Then    
$$H^+_0(  \underline S^{(0)}   | \underline T \setminus 
 \underline   S^{(0)}  )\ge \frac \zeta 2  |\G_0|^\a,  \Eq (fo.1) $$
where }
$$  |\G_0|^\a:=  \sum_{ T \in \underline  \G_0  }  |T|^\a.   \Eq (fo.5)  $$ 
\vskip0.5cm
 
  \noindent  { \bf \Theorem  (2CFMP)  [\rcite {CFMP}]} {\it  For  any    $\g>0$ 
  there exists  $C_0(\g)$ so that for $b\ge  C_0(\g)$  and for all $m>0$    
  
     $$ \sum_{0 \in \Gamma |\G| =m }   w_{b}^\g  (\Gamma)  \le 2m e^{- b m^{\g}},\Eq (E.6a) $$ 
where }
$$ w_b^\g (\Gamma):= \prod_{ T \in \G} e^{-b| T|^{\g}}.  \Eq (E.4a) $$ 
  \vskip0.5cm

  \noindent  In the sequel,  it is convenient   to identify  in    each contour $\G$  the  
families  of triangles having the same mass.   
  \vskip0.5cm
   \noindent {\bf \Definition  (2a)}
$$ \Gamma  =    \{\underline T^{(0)}, \underline T^{(1)},\dots \underline T^{(k_{\G})}  \}$$  
where   for    $\ell=0, \dots k_{\G}$,  $\underline T^{(\ell)}:= \{ T^{(\ell)}_1, T^{(\ell)}_2, 
\dots  T^{(\ell)}_{n_\ell} \}$,  and  each triangle of the family $\underline T^{(\ell)}$  has
 the same mass, i.e.   for all $i \in \{1, \dots n_\ell\} $, $| T^{(\ell)}_i| = \Delta_\ell$  for    
  $\Delta_\ell \in \N$.     According to  \eqv (fo.5) 
$$ |\G |^\r = \sum_{\ell=0}^{k_{\G }} |\underline T^{(\ell)} |^\r, \quad 
  |\underline T^{(\ell)} |^\r =\sum_{T\in \underline T^{(\ell)}} |T|^\r =  n_\ell \Delta_\ell^\r, 
\quad \r \in \R^+. \Eq (Ma.3) $$ 
\vskip 1truecm
\chap{3  Proof of Theorem  \eqv(1)  }3
\numsec= 3
\numfor= 1
\numtheo=1
\medskip \noindent 

The proof of  Theorem  \eqv(1)   is an immediate  consequence  of the following proposition  and the Markov inequality. 
\vskip0.5cm \noindent  
{\bf \Proposition (2)  } {\it     Let $ \alpha \in (\frac 12, \frac {\ln3} {\ln2} -1)$. 
There exist    positive  %$\zeta:= \zeta(\a)$,  see \eqv(fo.2),   
$\theta_0:= \theta_0(\alpha)>0 $ and  $ \b_0:=\b_0(\a)>0 $  so that for  
 $0<\theta \le \theta_0 $  and    $\b \ge  \b_0 $   
$$
\E\left [ \mu^+ \Big( \frac {1- \s_0} 2\Big)\right] 
\le   e^{-\frac {\bar b} {100}  }  , \Eq (R2)  
$$
where $\bar b$  is the quantity defined in \eqv (SS5). 
 } 

\medskip
\proof
A  necessary condition to have $\s_0=-1$ is that the site zero is contained in the support of some 
 contour $\G$ so that 
$$  \mu^+_{ \Lambda}(\s_0=-1) \le  \mu^+_{ \Lambda}(\{\exists \G : 0 \in \Gamma\}) \le  \sum_{\Gamma \ni 0 } 
 \mu^+_{ \Lambda}( \Gamma).  
\Eq(Basic1)
$$
  By definition, see \eqv (2.3),
 $$\mu^+_{ \Lambda}( \Gamma)[\o]:=  \frac 1 { Z^+_{\Lambda}[\o]}  \sum_{\underline T: 
\underline T \simeq \G    } e^{-\beta H^+(\underline T \cup   \G  )[\o]} ,   \Eq (E.8a) $$ 
where  $ \sum_{\underline T:  \underline T \simeq \G    }$ means that the  sum is %done 
over all    families  of triangles  compatible with the contour $\G$.
 Recalling \eqv (2.1a) and \eqv (D3), %picking up a
  for any   $j$ such that 
 $0 \le  j \le  k_{\G}$,  
 % that  $$  H^+(\underline T \cup     \G ):= H_{0}^{+ }(\underline T  \cup    \G)+  \theta  G 
%(\s(\underline T  \cup    \G)) . $$ 
 we write  for the deterministic part of the Hamiltonian  
  $$   
H_{0}^{+ }(\underline T  \cup   \G )=    
 H^+_{0 } (\underline T  \cup   \G \setminus ( \cup_{\ell=0}^j\underline T^{(\ell)})
  ) + H^+_0(\underline T \cup   \G \setminus ( \cup_{\ell=0}^j\underline T^{(\ell)})  |
  ( \cup_{\ell=0}^j\underline T^{(\ell)}) ).  
\Eq(I1)
$$ 
 Using  estimate   \eqv(fo.1) and recalling  notation  \eqv(Ma.3)   
$$
H^+_0(\underline T \cup   \G \setminus ( \cup_{\ell=0}^j\underline T^{(\ell)})  | 
( \cup_{\ell=0}^j\underline T^{(\ell)}) )  \ge\frac  \zeta 2  \sum_{\ell=0}^j  n_{\ell} 
|\Delta_{\ell}|^\alpha. 
\Eq (M1a)
$$
Therefore 
 $$   \mu^+_{ \Lambda}( \Gamma)      \le e^{-  \beta  \frac \zeta 2 \sum_{\ell=0}^j  n_{\ell} 
 |\Delta_{\ell}|^\alpha }     \frac 1 { Z^+_{\Lambda}}    \sum_{\underline T: \underline
 T \simeq \G }    e^{-\beta H^+_{0} (\underline T \cup   \G  \setminus ( \cup_{\ell=0}^j
\underline T^{(\ell)})) + \beta \theta G(\s (\underline T\cup    \G  ) )[\o]}. 
  \Eq (M.2) $$
 We multiply and divide \eqv (M.2)  
  % Multiplying  and dividing  \eqv (M.2) 
   by
 $$ \sum_{\underline T:\underline T \simeq \G  } e^{-\beta H^+_{0} (\underline T \cup   
 \G \setminus   (\cup_{\ell=0}^{j}  \underline T^{(\ell)})) + \beta \theta G(\s (\underline T
 \cup   \G  \setminus   \cup_{\ell=0}^j  \underline T^{(\ell)} ))[\o]  },    \Eq (L.1ba) $$ 
and   reconstruct  
 $\mu^+_{ \Lambda}( \Gamma\setminus  [ \cup_{\ell=0}^j \underline T^{(\ell)} ] )$,
 observing  that 
  $\sum_{\underline T: \underline T \simeq \G } 1   \le 
   \sum_{\underline T: \underline T \simeq \Gamma \setminus    \cup_{\ell=0}^j 
 \underline T^{(\ell)}  }  1   % \Eq (L.5)
  $.  
%one   
% so that 
  We get   $$
 \mu^+_{ \Lambda}( \Gamma) \le 
  e^{-\frac {\beta \zeta}2 \sum_{\ell=0}^j |\underline T^{(\ell)}|^\alpha }  
\mu^+_{ \Lambda}( \Gamma\setminus  [ \cup_{\ell=0}^j \underline T^{(\ell)} ] ) 
e^{\b  F_j[\o]} 
 \Eq (E.2a) $$ 
where 
%we  set 
  $$  F_j[\o]: =  \frac 1 \b \ln \left \{  \frac {   \sum_{\underline T: \underline T \simeq \G   } 
 e^{-\beta H^+_{0} (\underline T \cup  \G   \setminus    \cup_{\ell=0}^j  \underline T^{(\ell)})  +
 \beta \theta G(\s (\underline T  \cup   \G ) )[\o]}}
  {  \sum_{\underline T: \underline T \simeq \G   }  e^{-\beta H^+_{0} (\underline T \cup   
 \G \setminus   (\cup_{\ell=0}^j \underline T^{(\ell)} )) + \beta \theta  G 
 (\s (\underline T  \cup   \G \setminus   (\cup_{\ell=0}^j \underline T^{(\ell)} ))[\o]} }\right\} . 
 \Eq (L.1b) $$ 
   In  \eqv (E.2a) we explicitly quantify  the deterministic cost
of the   first smaller  families  of triangles   $\{ \underline T^{(0)}, \dots, \underline T^{(j)} \}$
and    express  the main random contribution  $  F_j[\o]$ so that it is antisymmetric  with respect   to the     sign  exchange of  the random field inside  $\cup_{\ell=0}^j \underline T^{(\ell)} $, see \eqv (Anti1). This   observation    allows to  estimate  this random contribution in a very convenient way, see Lemma \eqv(e2).   To this aim  
we define  for each $\G$  the partition:
$ \O= \cup_{j=-1}^{k_{\G}}  B_j$  
where 
 for $j\in \{0, \dots k_{\G}-1\}$ 
  $$B_j=B_j (\G):=\{ \o :   F_j[\o]  \le  \frac \zeta 4    
 \sum_{\ell=0}^j |\underline T^{(\ell)}|^{ \alpha},     \hbox { \rm and } \; 
\forall i >j,   F_i [\o]  >   \frac \zeta 4     \sum_{\ell=0}^i |\underline T^{(\ell)}|^{ \alpha}  \},  
 \Eq (P.1) $$ 
  $$B_ {k_{\G}}= B_ {k_{\G}} (\G) :=\{ \o :  F_{k_{\G}}[\o]  \le  \frac \zeta 4 
  \sum_{\ell=0}^{k_{\G}} |\underline T^{(\ell)}|^{ \alpha}  \},   \Eq (P.1a) $$ and 
   $$B_{-1} =  B_{-1} (\G) :=\{ \o :  \forall i >-1 ;   F_i[\o]     >   \frac \zeta 4  
 \sum_{\ell=0}^i |\underline T^{(\ell)}|^{ \alpha}  \}.  \Eq (P.2) $$ 
 The  relevant properties  of the partition are given in  the following lemma,  whose  proof is given 
  in    Section 4.  
\medskip
\noindent{\bf \Lemma(e2)}{\it 
For  $-1 \le  j \le  k_{\G}$,   
 $$ \E \left [ \1_{B_j} \right ] \le e^{-   \frac {\z^2 }{2^{10}\theta^2}   
  \sum_{\ell =j+1}^{ k_{\G}} |\underline T^{(\ell)}|^{2\alpha-1} }   .  \Eq (E.3a)$$
with the convention that an empty sum is zero.}

\medskip
\noindent
  We then write 
  $$ \mu^+_{ \Lambda}( \Gamma)= \sum_{j=-1}^{ k_{\G}}  \mu^+_{ \Lambda}( \Gamma)\1_{\{B_j\}}
   $$ 
and  apply to each   $ \mu^+_{ \Lambda}( \Gamma)\1_{\{B_j\}}$
estimate  \eqv (E.2a).   
 %using Lemma \eqv (e2),
   We obtain
$$
\eqalign { 
  \E \left  [\mu^+_{ \Lambda}( \Gamma)\right] &=   \sum_{j=-1}^{ k_{\G}}   
\E \left  [\mu^+_{ \Lambda}( \Gamma) )\1_{\{B_j\}} \right]  \cr & \le \sum_{j=-1}^{ k_{\G}} 
   e^{-\frac {\beta \zeta} 4 \sum_{\ell=0}^j |\underline T^{(\ell)}|^\alpha }  
e^{- c \frac {\z^2 }{\theta^2}  \sum_{\ell =j+1}^{ k_{\G}} |\underline T^{(\ell)}|^{2\alpha-1} } \cr &\le
(k_\G +1) e^{-\bar b  \sum_{\ell =0}^{ k_{\G}} |\underline T^{(\ell)}|^{2\alpha-1} },
} \Eq   (En.1) $$
where 
$ \bar b:= \min(\frac {\beta \zeta} 4, \frac {\z^2} {2^{10}\theta^2}).  $
Recalling \eqv (E.4a), one has
$$ 
\E \left [  \mu^+_{ \Lambda}(\{ 0 \in \Gamma\})\right]  \le    \sum_{ \Gamma \ni 0 }  
(k_{\G}+1)  w_{\bar b}^{2\a-1}  (\Gamma)=   \sum_{m \ge 3}  (m+1) 
 \sum_{0 \in \Gamma: |\G| = m }   w_{\bar b}^{2\a-1}
  (\Gamma).
\Eq(est)
$$
Using \eqv(E.6a), after  a few lines computation one gets \eqv(R2)\eop
\medskip
    \noindent {\bf Remark:}   The upper bound $ \a < \frac {\ln 3 } {\ln 2} -1$  in Theorem \eqv (1) follows from Theorem  \eqv  (1CFMP), the lower bound $ \a >\frac 12$ from Theorem \eqv  (2CFMP)
and \eqv  (est).   
\vskip 1truecm
\chap{4  Probabilistic estimates
}4
\numsec= 4
\numfor= 1
\numtheo=1 
 \noindent 

Let   $h=h[\o]$ be a  realization of the random magnetic fields and 
$A\subset \Z$.  Define   
$$
(S_{A}h)_i=
\cases {-h_i,&if $i\in A$;\cr
h_i, &otherwise ,\cr}
\Eq(sf)
$$
and denote $h[S_A \o]\equiv S_Ah[\o]$.    In the following to simplify notation we wet
$  S_{\underline T} h= S_{{\rm supp} (T)}h $.
Recalling \eqv(2.1ab), it is easy to see that
$$ 
 G(\s (\underline T  \cup  \G  \setminus   \underline T^{(0)}  ))[\o] 
=G(\s(\underline T\cup \G))[S_{\underline T^{(0)}}\o].
\Eq(basic1)
$$
In general 
$$   
G(\s 
(\underline T   \cup  \G  \setminus \cup_{\ell=0}^i  \underline T^{(\ell)}  ))[\o]= 
G(\s (\underline T  \cup  \G)) [ S_{D_i} \o] 
\Eq(basic2)
$$ 
where 
$$
D_i\subset \cup_{\ell=0}^i \left ( {\rm supp} (\underline T^{(\ell)})\right )
\Eq(d1)
$$   
 is the non--empty set so  that
$$
S_{D_i}=S_{\underline T^{(i)}} S_{\underline T^{(i-1)}} \dots 
S_{\underline T^{(1)}} S_{\underline T^{(0)}}. 
\Eq(basic3)
$$
When all the triangles in  $ (\underline T^{(\ell)}, \ell=0,\dots,j)$ have  disjoint supports 
 \eqv(d1) becomes  an equality. 
In general there are triangles inside triangles and  in this case 
the inclusion in \eqv(d1) is strict.
 By  construction  
the     $  F_j[\o] $ defined in    \eqv (L.1b)    are  such that 
  $$ F_j (h(D_j^c), h(D_j))= - F_j (h(D_j^c), -h(D_j)),  \qquad  j \in 0, \dots, k_\G, \Eq (Anti1)$$
where for a set $A  \subset  \Z$, we denote by $ h(A) = \{ h_i: i \in  A\}$.
Therefore one gets   that $\E [ F_j]=0$.

        \medskip
  \noindent{\bf Proof of Lemma \eqv (e2) }   Set 
$$  A_i:=   \frac \zeta 4    \sum_{\ell=0}^i |\underline T^{(\ell)}|^{ \alpha}.  \Eq (T1) $$
We have 
$\P \left [ B_j  \right ]    \le  \P \left [  \forall i >j;     F_i[\o]     >   A_i    \right ]. $
  Let $ \l_i $ for  $i=j+1, \dots,  k_\G$ be      positive parameters,  
 by exponential Markov inequality we have 
 $$
   \P \left [ \forall i >j: F_i [\o] \ge  A_i    \right] 
\le  e^{-   \sum_{\ell=j+1}^{k_\G}   \l_\ell  A_\ell }  
 \E\left [ e^{ \sum_{\ell=j+1}^{k_\G} \l_{\ell}   F_\ell  }    \right ] .
 \Eq (V.5a) 
$$ 
Set
$$ 
F[\o] :=  \sum_{i=j+1}^{k_\G} \l_{i}  F_i[\o]. 
\Eq(effe)
$$
 It remains to estimate $\E[e^F]$.
Note that $ F[\o]$ depends on all  the random fields on $\Lambda$. Let $N$ be the number of sites in $\L$.
To avoid involved notations, we define a bijection $\Pi$ from $\L$ to $\{1,\dots,N\}$ as follows:
first pick up all the $n_0\D_0$ sites in  $\hbox {supp}(\underline T^{(0)})$
and put them consecutively in $N,\dots, N-n_0\D_0+1$ 
(keeping them in the same order as they are for definiteness). Then pick up the sites 
in  $\hbox {supp}(\underline T^{(1)})$ that are   not  in  $\hbox {supp}(\underline T^{(0)})$
and put  them consecutively starting at $N-n_0\D_0$ until they are exhausted.
Note that if no triangles of size $\D_0$ are within triangle of size $\D_1$, $\Pi$ maps
${supp}(\underline T^{(0)})\cup {supp}(\underline T^{(1)})$ onto $\{N,\dots,N-n_0\D_0-n_1\D_1+1\}$;
otherwise $\Pi$ maps 
${supp}(\underline T^{(0)})\cup {supp}(\underline T^{(1)})$ onto a proper subset of 
 $\{N,\dots,N-n_0\D_0-n_1\D_1+1\}$.
One can iterate this procedure until all the sites of the support of $\G$ are  exhausted.
As above, for all $1\le j\le k_\G-1$,  if all triangles considered are disjoint $\Pi$ maps 
$\cup_{\ell=0}^{j+1} {supp}(\underline T^{(\ell)})$ onto
$\{N,N-1,\dots, N-M_{j+1}+1\}$ where 
$ M_{j+1}=   \sum_{\ell=0}^{j+1} n_\ell \D_\ell$,   
 otherwise on a proper subset of it.
Then one can pick up all the remaining sites of $\L$ and continue as above.   The $\Pi$ so defined   induces a bijection from  the random magnetic fields indexed by $\L$ to a family of random variables 
$(h_1,\dots,h_N)$ by  $(\Pi h)_i:=h_{\Pi i}, \forall i\in \L$.
Using this bijection, one can work with the random variables  $(h_i, 1\le i\le N)$.
Define  the  family of increasing $\s$-algebra:
$$
\left(\emptyset,\O\right)=\S_0\subset\S_1=\s(h_1)\subset\S_2=\s(h_1,h_2)\subset   
\dots \subset \S_N=\s(h_1,h_2,\dots,h_N) $$
and 
$
\D_k(F)=\E\left[F|\S_k\right]-\E\left[F|\S_{k-1}\right]
 $ the associated martingale difference sequences.
We have 
 $$ \E [ F| \S_N] = F;\quad   \E [ F|\S_0] = \E [F]=0, \quad F= \sum_{k=1}^N \D_k(F). $$
 Remark  that 
 $$ 
\E  [ F_{j+1} |\S_{i}] = 0  \quad  \forall  i \in \{1,\dots  N- M_{j+1} \}
\Eq(fonda1)   
$$
 since by \eqv(Anti1) 
$$
F_{j+1}(h(D^c_{j+1}),h(D_{j+1}))=-F_{j+1}(h(D^c_{j+1}),-h(D_{j+1}))
$$
and 
$$\Pi h(\cup_{\ell=0}^{j+1} \underline T^{(\ell)} ) \subset 
\{ h_{N-M_{j+1}},\dots, h_N \}.$$
Note also that 
 $$ 
\E  [ F |\S_{i}] = 0  \quad  \forall  i \in \{1,\dots  N- M_{k_\G} \},  
\Eq(fond1)
$$
and 
$$ 
\E[e^{F}]=\E[e^{\sum_{k=1}^{N-1}\D_k(F)} \E[e^{\D_N(F)}|\S_{N-1}]].
$$
With self explained notations, using Jensen inequality one has
$$
\E[e^{\D_N(F)}|\S_{N-1}]=
\int e^{ \D_N(F)}\,\P(dh_N)\le 
\int e^{ \left[F(h_{<N},h_N)-F(h_{<N},\tilde h_N)\right]} \P(dh_N)\P(d\tilde h_N).
\Eq(4.20a)
$$
We then expand   the exponential in the right hand side of \eqv(4.20a).    All the odd powers but 
the constant one vanish.   For the even power  we recall    \eqv(effe) and
by the  Lipschitz  continuity of each term  with respect to $(h_1,\dots, h_N)$ we get 
  $$\eqalign {    \left  | F(h_{<N},h_N)-F(h_{<N},\tilde h_N)\right| & \le     \sum_{i=j+1}^{k_\G} \l_{i}  | F_i (h_{<N},h_N) -  F_i (h_{<N},\tilde h_N)|\cr &  \le 
   2\th |h_N-\tilde h_N|
\sum_{\ell=j+1}^{k_\G} \l_{\ell}. }
\Eq (A.5)
 $$  
Then  estimating  $|h_N-\tilde h_N|\le 2$
and  $2^{(2n-1)^+}(2n!)^{-1}\le (n!)^{-1}$ to re-sum the series   one gets
$$
\E[e^{\D_N}|\S_{N-1}]\le e^{16\th^2 ( \sum_{\ell=j+1}^{k_\G} \l_{\ell})^2}.
\Eq(step1)
$$
In the case of gaussian or  subgaussian variable one just performs all 
the integration instead of using $|h_N-\tilde h_N|\le 2$. This will  modify the result by a constant 
different from 16.
To iterate,  one uses again the Jensen inequality  obtaining 
$$
\eqalign{
\E[e^{\D_{N-1}}|\S_{N-2}]
&=\int e^{ \D_{N-1}(F)}\P(dh_{N-1})\cr
&\le 
\int e^{
\int \left[F(h_{<N-1},h_{N-1},\hat h_N)-
F(h_{<N-1},\tilde h_{N-1},\hat h_N)\right]\P(d\hat h_N)
}
\P(dh_{N-1})\P(d\tilde h_{N-1}). 
} 
$$
It is clear that the random variable
$$
\int \left[F(h_{<N-1},h_{N-1},\hat h_N)-
F(h_{<N-1},\tilde h_{N-1},\hat h_N)\right]\P(d\hat h_N)
$$
is a symmetric ones under $\P(dh_{N-1})\P(d\tilde h_{N-1})$ and satisfies an estimate as \eqv(A.5)
from which one gets
$$
\E[e^{\D_{N-1}(F)}|\S_{N-2}]\le e^{16\th^2 ( \sum_{\ell=j+1}^{k_\G} \l_{\ell})^2}.
\Eq(step2)
$$
Iterating one gets  $\E[e^{\D_{k}(F)}|\S_{k-1}]\le 
e^{4\th^2 ( \sum_{\ell=j+1}^{k_\G} \l_{\ell})^2}$ for 
$k\in \{N,N-1,\dots, N-M_{j+1}\}$. 
When $k=N-M_{j+1}-1$, a new fact happens.   Using \eqv(fonda1) for $i=N-M_{j+1}\equiv m $ 
and computing 
$$
\D_m(F)=\sum_{i=j+1}^{k_\G} \l_{i} \left(\E[F_i|\S_m]-\E[F_i|\S_{m-1}]\right)
\Eq(fond2)
$$
one  obtains  that the term corresponding to $i=j+1$ in the sum gives   zero contribution. 
Therefore, in this case, one  has
$$\left |  \int \left[F(h_{<m},h_m,  \hat  h_{>m})-
F(h_{<m}, \tilde h_m,  \hat  h_{>m})\right]  \P(\hat  h_{>m}) \right |   \le   
 4  \theta \sum_{\ell=j+2}^{k_\G} \l_{\ell}.   
\Eq(fond3)
$$
Iterating  this procedure one gets
$$
\E\left [ e^{ F }    \right ] \le  e^{ 16   \theta^2 \left 
\{ \left ( \sum_{\ell=0}^{j+1} n_{\ell}\Delta_{\ell} \right ) 
\left (  \sum_{\ell=j+1}^{k_\G} \l_{\ell}\right )^2     + 
 n_{j+2}\Delta_{j+2} \left (  \sum_{\ell=j +2}^{k_\G} \l_{\ell}\right )^2+  
 \dots +  n_{k_\G}\Delta_{k_\G} \left (    \l_{k_\G}\right )^2\right \}  } . 
\Eq (V.5b)  
$$
The estimate \eqv(V.5b) suggests to set for $\ell=j+1,\dots,k_\G$
$$
\mu_\ell\equiv    \sum_{n= \ell}^{k_\G}  \l_{n}
\Eq (Pa.1) 
$$ 
and the constraints $(\l_i\ge 0, j+1\le i\le k_\G)$ become $\mu_\ell$ decreasing with $\ell$.
We write  the first exponent of  \eqv(V.5a)    in terms of $\{\mu_\ell\}_{\ell=0}^j$ obtaining  
$$ 
-\sum_{\ell=j+1}^{k_\G}   \l_\ell  A_\ell =  
 - \frac \zeta 4    
\mu_{j+1}   \left (\sum_{\ell=0}^{j} n_{\ell}\Delta^{\a}_{\ell} \right ) -
\sum_{\ell=j+1}^{k_\G}   \frac \zeta 4    
\mu_{\ell} n_{\ell} \Delta^{\a}_{\ell}, 
\Eq(new1)
$$ 
and for  the exponent in \eqv(V.5b)  we obtain 
$$
16 \th^2 (\mu_{j+1})^2\left(\sum_{\ell=0}^j n_\ell\D_\ell\right)
+
\sum_{\ell=j+1}^{k_\G} 16 \th^2 (\mu_{\ell})^2 n_{\ell} \Delta_{\ell}.
\Eq(new2)
$$
Denote 
 $$
 f(\mu_\ell) \equiv -  \frac \zeta 4     \mu_\ell     \Delta_{\ell}^{ \alpha}    + 16  \theta^2 
   \mu^2_\ell   \Delta_{\ell}   \quad \ell = j+1,\dots  k_\G . $$ 
Choose  $\mu_\ell\equiv \bar \mu_{\ell}$ where 
 $$\bar \mu_{\ell} =   \frac 1 {4 \times 32} \frac \zeta   {  \th^2  \Delta_\ell^{1-\alpha } },  
\Eq (final1) $$
is the minimizer of  $f(\mu_\ell)$. Note that  $\bar \mu_{\ell}$
   is a  decreasing  function of  $\ell$    
and  
$$
f(\bar \mu_{\ell})=-\frac{\z^2\D_\ell^{2\a-1}}{2^{10}\th^2}.
\Eq(prefinal)
$$
Therefore  collecting together the  last sum in \eqv(new1) and the one  in \eqv(new2)    we get 
$$
-\sum_{\ell=j+1}^{k_\G}   \frac \zeta 4    
\bar\mu_{\ell} n_{\ell} \Delta^{\a}_{\ell}+
\sum_{\ell=j+1}^{k_\G} 16 \th^2 (\bar \mu_{\ell})^2 n_{\ell} \Delta_{\ell}
= -\frac{\z^2}{2^{10}\th^2} \sum_{\ell=j+1}^{k_\G} n_{\ell}\Delta^{2\a-1}_\ell 
=-\frac{\z^2}{2^{10}\th^2}\sum_{\ell=j+1}^{k_\G} |\underline T^{(\ell)}|^{2\a-1}.
\Eq(prefinal2)
$$
Summing up   \eqv(new1) and  \eqv(new2), taking in account \eqv (final1) and \eqv (prefinal2)  we get 
 $$
%\eqalign{
%&
-\frac \zeta 4    
\bar \mu_{j+1}   \left (\sum_{\ell=0}^{j} n_{\ell}\Delta^{\a}_{\ell} \right ) +
16 \th^2 (\bar \mu_{j+1})^2\left(\sum_{\ell=0}^j n_\ell\D_\ell\right)
%\cr &
=-\sum_{\ell=0}^jn_{\ell}
\left(\frac{\z}{4} \bar \mu_{j+1} \Delta^\a_\ell-16\th^2\Delta_\ell (\bar\mu_{j+1})^2\right).
%\cr }
\Eq(prefinal3)
$$
One can check easily that for all $0\le \ell \le j $ one has  
$$
\left(\frac{\z}{4} \bar \mu_{j+1} \Delta^\a_\ell-16\th^2\Delta_\ell (\bar\mu_{j+1})^2\right)
=\frac{\z^2\Delta_{\ell}}{2^9\th^2\Delta_{j+1}^{1-\a}}
\left(\frac 1{\Delta_\ell^{1-\a}}-\frac1{2\Delta_{j+1}^{1-\a}}\right)
>0
\Eq(prefinal4)
$$
since by construction $\Delta_\ell<\Delta_{j+1}$ for $0\le \ell\le j$. 
%Putting together
  By   \eqv(V.5a), \eqv(V.5b), \eqv(prefinal2),  and \eqv(prefinal3) one gets \eqv(E.3a).\eop

\vskip 0.5cm
 \noindent{\bf Acknowledgements}
We  are indebted to  Errico Presutti for stimulating discussions and criticism. 
P.P.  thanks   the   Mathematics  Department   of `` Universit\'a degli Studi dell'Aquila'' 
and  Anna de Masi for hospitality. 
Enza Orlandi   thanks the   Institut Henri
Poincar\'e - Centre Emile Borel, (workshop 
M\'ecanique statistique, probabilit\'es et syst\`emes de particules 2008)
for hospitality.  
The authors thank the referees for useful comments.

\centerline{\bf References}
\vskip.3truecm
 
\item{[\rtag{ACCN}]} M. Aizenman, J. Chayes, L. Chayes and C. Newman:
Discontinuity of the magnetization in one--dimensio\-nal $1/|x-y|^2$ percolation, Ising and Potts models.
{\it J. Stat. Phys.} {\bf 50} no. 1-2 1--40 (1988).

 \item{[\rtag{AW}]} M. Aizenman, and  J. Wehr:
{ Rounding of first or\-der pha\-se tran\-si\-tions
in sys\-tems with quenched disorder.}
{\it Comm. Math. Phys.} {\bf 130}, 489--528 (1990).

 \item{[\rtag{Bo}]} A. Bovier: 
{\it  Statistical Mechanics of Disordered Systems.}
{Cambridge Series in Statistical and Probabilistic mathematics}, (2006).
 
 \item{[\rtag{BK}]} J. Bricmont, and A. Kupiainen:
{  Phase transition in the three-dimensional
random field Ising model.}
{\it Comm. Math. Phys.},{\bf 116}, 539--572 (1988).

 \item{[\rtag {CFMP}]} M. Cassandro, P.  A. Ferrari, I. Merola and E. Presutti:  
  { Geometry of contours and Peierls estimates in $d=1$ Ising models with long range 
interaction.} {\it J. Math. Phys.} {\bf 46},   no 5,  053305,   (2005).

\item{[\rtag{D0}]} R. Dobrushin: 
The description of a random field by means of conditional probabilities and. conditions of its regularity.
{\it  Theory Probability Appl.} {\bf 13}, 197-224 (1968)  
\item{[\rtag{D1}]} R. Dobrushin:
The conditions of absence of phase transitions in one-dimensional classical systems.
{\it Matem. Sbornik}, {\bf 93} (1974), N1, 29-49

\item{[\rtag{D2}]} R. Dobrushin:
Analyticity of correlation functions in one-dimensional classical systems with slowly decreasing potentials.
{\it Comm. Math. Phys.} {\bf 32} (1973), N4, 269-289

  \item{[\rtag{Dy}]} F.J. Dyson:
 {  Existence of  phase transition in a  one-dimensional Ising ferromagnetic.} 
{ \it Comm. Math. Phys.},{\bf 12},91--107,  (1969).

\item{[\rtag{E}]} R.S. Ellis: {\it Entropy, Large deviation and Statistical mechanics.}
 New York: Springer (1988).

\item{[\rtag{FS}]} J.  Fr\"ohlich and  T. Spencer:  
 {  The phase transition in the one-dimensional Ising model with
 $\frac 1 {r^2}$ interaction energy.} { \it Comm. Math. Phys.}, {\bf 84}, 87--101,   (1982). 
  
\item{[\rtag{GM}]} G. Gallavotti and S. Miracle Sole:
Statistical mechanics of Lattice Systems. 
{\it Comm. Math. Phys.} {\bf 5} 317--323 (1967)

 \item{[\rtag{IM}]} Y. Imry and  S. Ma:  
 {  Random field instability of the ordered state of continuous symmetry.} 
{ \it  Phys. Rev. Lett.},  {\bf 35}, 1399--1401,   (1975).

\item {[\rtag{RT}]} J. B.  Rogers and C.J. Thompson:  
Absence of long range order in one dimensional spin systems. 
{\it  J. Statist. Phys.} {\bf 25}, 669--678 (1981)

\item {[\rtag{Ru}]} D. Ruelle:  
Statistical mechanics of one-dimensional Lattice gas.
{\it  Comm. Math. Phys.} {\bf 9}, 267--278 (1968)

   \end